\documentclass{gtart}

\def\ifplaintex{\expandafter\ifx\csname documentclass\endcsname\relax}

\ifplaintex 
\else
\fi

\def\gtm{{\mathsurround=0pt\it $\cal G\mskip-2mu$eometry \&\ 
$\cal T\!\!$opology $\cal M\mskip-1mu$onographs}}    

\def\gtp{{\mathsurround=0pt\it $\cal G\mskip-2mu$eometry \&\ 
$\cal T\!\!$opology $\cal P\!$ublications}}  

\def\recd{{\small Received:\qua\receiveddate\ifx\reviseddate\relax
\else\qquad Revised:\qua\reviseddate\fi\par}} 

\def\volumenumber#1{\def\thevolumenumber{#1}}
\def\volumeyear#1{\def\thevolumeyear{#1}}
\def\volumename#1{\def\thevolumename{#1}}
\def\papernumber#1{\def\thepapernumber{#1}}
\def\pagenumbers#1#2{\def\startpage{#1}\def\finishpage{#2}}
\def\published#1{\def\publishdate{#1}}
\def\received#1{\def\receiveddate{#1}}
\def\revised#1{\def\reviseddate{#1}}
\def\accepted#1{\def\accepteddate{#1}}
\def\asciititle#1{\def\theasciititle{#1}}

\def\coverauthors#1{\def\thecoverauthors{#1}}
\def\asciiauthors#1{\def\theasciiauthors{#1}}
\def\asciiaddress#1{\def\theasciiaddress{#1}}

\def\coverauthors#1{\def\thecoverauthors{#1}}
\long\def\asciiabstract#1{\long\def\theasciiabstract{#1}}
\def\asciikeywords#1{\def\theasciikeywords{#1}}

\let\\\par
\let\thevolumenumber\relax\let\thepapernumber\relax
\let\thevolumeyear\relax\let\startpage\relax
\let\finishpage\relax\let\publishdate\relax\let\receiveddate\relax
\let\reviseddate\relax\let\accepteddate\relax\let\theasciititle\relax
\let\theasciiauthors\relax\let\theasciiaddress\relax
\let\theasciiabstract\relax\let\theasciikeywords\relax
\let\thecoverauthors\relax
\let\thecoverauthors\relax\let\theerratum\relax\let\theasciiemail\relax
\let\theshortauthors\relax\let\theshorttitle\relax

\def\startpage{1}\def\finishpage{15}\def\thepapernumber{77}

\volumenumber{2}
\volumename{Proceedings of the Kirbyfest}
\volumeyear{1999}

\long\def\maketitlep{   

\count0=\startpage

\gtm\nl        
{\small Volume \thevolumenumber: \thevolumename\nl 
\ifx\theerratum\relax\else Erratum \erratumnumber\nl\fi
Pages \startpage--\finishpage\nl}

\vglue 0.1truein   

{\parskip=0pt\leftskip 0pt plus 1fil\def\\{\par\smallskip}{\ifplaintex\large
\else\Large\fi\bf\thetitle}\par\medskip}   
\vglue 0.05truein 

{\parskip=0pt\leftskip 0pt plus 1fil\def\\{\par}{\sc\theauthors}
\par\medskip}%
 
\vglue 0.03truein 

{\small\leftskip 25pt\rightskip 25pt{\bf Abstract}\stdspace\theabstract

{\bf AMS Classification}\stdspace\theprimaryclass
\ifx\thesecondaryclass\relax\else; \thesecondaryclass\fi\par
{\bf Keywords}\stdspace \thekeywords\par}\vglue 7pt

}   

\font\phead=cmsl9 scaled 950
\font=cmsl9 scaled 1050
\font\pnum=cmbx10 scaled 913
\font\lnum=cmbx10 
\font\pfoot=cmsl9 scaled 950
\font=cmsl9 scaled 1050
\ifplaintex
\headline{\vbox to 0pt{\vskip -4.5mm\line{\small\phead\ifnum
\count0=\startpage ISSN 1464-8997 (on line)
1464-8989 (printed) \hfill {\pnum\folio}\else\ifodd\count0\def\\{ }%
\ifx\theshorttitle\relax\thetitle\else\theshorttitle\fi\hfill{\pnum\folio}
\else\def\\{ and }{\pnum\folio}\hfill\ifx\theshortauthors\relax\theauthors
\else\theshortauthors\fi\fi\fi}\vss}}
\footline{\vbox to 0pt{\vglue 0mm\line{\small\pfoot\ifnum\count0=\startpage
Published \publishdate:\qua\copyright\ \gtp\hfill\else
\gtm, Volume \thevolumenumber\ (\thevolumeyear)\hfill\fi}\vss
}}
\else
\makeatletter
\def\@oddhead{{\small\lhead\ifnum\count0=\startpage ISSN 1464-8997 (on line)
1464-8989 (printed) \hfill {\lnum\number\count0}\else\ifodd\count0
\def\\{ }\ifx\theshorttitle\relax \thetitle \else\theshorttitle\fi\hfill
{\lnum\number\count0}\else\def\\{ and }{\lnum\number\count0}
\hfill\ifx\theshortauthors\relax 
\theauthors\else\theshortauthors\fi\fi\fi}}\def\@evenhead{@oddhead}
\def\@oddfoot{\small\lfoot\ifnum\count0=\startpage Published \publishdate:\qua\copyright\ \gtp\hfill\else
\gtm, Volume \thevolumenumber\ (\thevolumeyear)\hfill\fi}
\def\@evenfoot{@oddfoot}
\makeatother
\fi

\let\maketitlepage\maketitlep
\let\makeshorttitle\maketitlepage
\let\maketitle\maketitlepage

\newwrite\gtoutfile
\long\gdef\makeheadfile{  
{\def\\{, }\def\s{ }
\immediate\openout\gtoutfile head.xxx
\immediate\write\gtoutfile{To: math@arxiv.org}
\immediate\write\gtoutfile{Subject: put OR rep NNNNN:ppppp}
\immediate\write\gtoutfile{--text follows this line--}
\immediate\write\gtoutfile{Proxy-for: \ifx\theasciiauthors\relax
\theauthors\else\theasciiauthors\fi\s<\ifx\theasciiemail\relax\theemail\else\theasciiemail\fi>}
\immediate\write\gtoutfile{\noexpand\\}
\immediate\write\gtoutfile{Authors: \ifx\theasciiauthors\relax
\theauthors\else\theasciiauthors\fi}
{\def\\{ }\immediate\write\gtoutfile{Title: \ifx\theasciititle\relax
\thetitle\else\theasciititle\fi}}
\immediate\write\gtoutfile{Subj-class: GT or SG, GR etc}
\immediate\write\gtoutfile{MSC-class: \theprimaryclass\ifx\thesecondaryclass\relax\else, \thesecondaryclass\fi}
\immediate\write\gtoutfile{Journal-ref: Geom. Topol. Monogr. \thevolumenumber\s
(\thevolumeyear) \startpage-\finishpage}
\immediate\write\gtoutfile{Comments: Published by Geometry and Topology Monographs at}
\immediate\write\gtoutfile{\s\s\s  http://www.maths.warwick.ac.uk/gt/GTMon\thevolumenumber/paper\thepapernumber.abs.html}
\immediate\write\gtoutfile{\noexpand\\}
\immediate\write\gtoutfile{}
\ifx\theasciiabstract\relax
\immediate\write\gtoutfile{\theabstract}\else
\immediate\write\gtoutfile{\theasciiabstract}\fi
\immediate\write\gtoutfile{}
\immediate\write\gtoutfile{\noexpand\\}
\immediate\write\gtoutfile{}
\immediate\closeout\gtoutfile}}  

\def\maketitlepage{\maketitlep\makeheadfile}
\let\makeshorttitle\maketitlepage
\let\maketitle\maketitlepage

\volumenumber{4}
\volumename{Invariants of knots and 3-manifolds (Kyoto 2001)}
\volumeyear{2002}
\papernumber{2}
\pagenumbers{13}{28}
\received{27 November 2001}
\revised{17 April 2002}
\accepted{22 July 2002}
\published{19 September 2002}

\usepackage{amssymb,amsmath,graphicx}

\newtheorem{teo}{Theorem}[section]

\newtheorem{lem-defi}[teo]{Lemma-Definition}

\newtheorem{conge}[teo]{Conjecture}

\newcommand{\mr}{\mathbb{R}}
\newcommand{\mc}{\mathbb{C}}
\newcommand{\mz}{\mathbb{Z}}
\newcommand{\mh}{\mathbb{H}}

\newcommand{\mn}{\mathbb{N}}

\newtheorem{theo}{Theorem}[section]

\newtheorem{proposi}[theo]{Proposition}

\newcommand{\Dd}{{\mathcal D}}
\newcommand{\Ee}{{\mathcal E}}
\newcommand{\Ff}{{\mathcal F}}

\newcommand{\Ii}{{\mathcal I}}

\newcommand{\Pp}{{\mathcal P}}

\newcommand{\Tt}{{\mathcal T}}

\newcommand{\Ww}{{\mathcal W}}
\newcommand{\Vv}{{\mathcal V}}

\begin{document}

\title{QHI, $3$-manifolds scissors congruence classes\\and the volume
conjecture}
\asciititle{QHI, 3-manifolds scissors congruence classes and the volume
conjecture}
\authors{St\'ephane Baseilhac\\Riccardo Benedetti}
\coverauthors{St\noexpand\'ephane Baseilhac\\Riccardo Benedetti}
\asciiauthors{Stephane Baseilhac, Riccardo Benedetti}

\address{Dipartimento di Matematica, Universit\`a di Pisa\\Via F. Buonarroti, 2, I-56127 Pisa, Italy}

\asciiaddress{Dipartimento di Matematica, Universita di Pisa\\Via F. Buonarroti, 2, I-56127 Pisa, Italy}

\email{baseilha@mail.dm.unipi.it, benedett@dm.unipi.it}

\begin{abstract} This is a survey of our work on 
Quantum Hyperbolic Invariants (QHI) of $3$-manifolds. We explain how
the theory of scissors congruence classes is a powerful geometric
framework for QHI and for a `Volume Conjecture' to make sense.
\end{abstract}
\asciiabstract{This is a survey of our work on 
Quantum Hyperbolic Invariants (QHI) of $3$-manifolds. We explain how
the theory of scissors congruence classes is a powerful geometric
framework for QHI and for a `Volume Conjecture' to make sense.}


\asciikeywords{Volume conjecture, hyperbolic 3-manifolds, scissors
congruence classes, state sum invariants, 6j-symbols, quantum
dilogarithm}

\maketitle

\section{Introduction}\label{intro}
This text is based on the talk that the second author gave at the
workshop {\it Invariants of Knots and $3$-Manifolds} (RIMS, Kyoto
2001, September 17 -- September 21), completed by some ``private'' talks
he gave at the same occasion. It reports on our joint works in progress. We
refer to \cite{BB1,BB0} for more details; in fact, in \cite{BB0} we develop
the ideas of sections 7--9 of \cite{BB1}, with some important differences
in the way they are concretized. Here we content ourselves with
providing precise definitions and statements. This is
summarized as follows.

\noindent Let $(W,L,\rho)$ be a triple formed by a smooth compact closed
oriented $3$-manifold $W$, a non-empty link $L$ in $W$ and a flat
principal $B$-bundle $\rho$ on $W$. We denote by $B$ the Borel
subgroup of upper triangular matrices of $SL(2,\mc)$.

\noindent One associates to $(W,L,\rho)$ a $\Dd$-{\it scissors
congruence class} $\mathfrak{c}_{\Dd}(W,L,\rho)$ which belongs to a
{\it (pre)-Bloch-like group} $\Pp(\Dd)$ built on suitably decorated
tetrahedra. The class $\mathfrak{c}_{\Dd}(W,L,\rho)$ may be
represented geometrically by any $\Dd$-{\it triangulation} of
$(W,L,\rho)$ and it depends on the topology and the geometry of the
triple $(W,L,\rho)$. For any odd integer $N>1$ and for any
$\Dd$-triangulation $\Tt$, one defines a {\it reduction} mod$(N)$
$\Tt_N$ of $\Tt$. It is obtained via a ``quantization'' procedure
using the cyclic representation theory of a quantum Borel subalgebra
$\Ww_N$ of $U_q (sl(2,\mc))$ specialized at the root of unit $\omega_N
= \exp (2\pi i/N)$. One of the main contributions of this work is to
have pointed out the relationship between this representation theory
and $B$-flat bundles, which are encoded in $\Tt$ by so called
simplicial ``full'' 1-cocycles - see Section \ref{dtet}. Basically,
this relationship relies on the theory of quantum coadjoint action of
\cite{DCP}, which should ultimately lead to generalizations of the
QHI, replacing $B$ by other algebraic Lie groups.

\noindent One also defines a family of complex valued {\it Quantum
Hyperbolic Invariants} (QHI) $K_N(W,L,\rho)$, which have  {\it state
sum} expressions $K(\Tt_N)$ based on any $\Tt_N$. The elementary
building blocks of $K(\Tt_N)$ are the {\it cyclic $6j$-symbols} of
$\Ww_N$. Roughly speaking $K_N(W,L,\rho)$ may be considered as a function of
the class $\mathfrak{c}_{\Dd}(W,L,\rho)$, which we shall write
$K_N(\mathfrak{c}_{\Dd}(W,L,\rho))$ (strictly speaking this is not
completely correct - see the end of \S \ref{QHI}). It turns out
that when $\rho$ is the trivial flat $B$-bundle one recovers the 
topological invariant conjectured in \cite{K1}.

\noindent The asymptotic behaviour of $K_N(W,L,\rho)$ when $N\to
\infty$ should depend on the $\Dd$-scissors congruence class
$\mathfrak{c}_{\Dd}(W,L,\rho)$.  In fact, to any triple $(W,L,\rho)$
one can associate also an $\Ii$-{\it scissors congruence class}
$\mathfrak{c}_{\Ii}(W,L,\rho)$ that belongs to $\Pp(\Ii)$, which is an
enriched version of the classical (pre)-Bloch group built on
hyperbolic ideal tetrahedra. The class $\mathfrak{c}_{\Ii}(W,L,\rho)$
may be represented by an $\Ii$-{triangulation} $\Tt_{\Ii}$ of
$(W,L,\rho)$, which is obtained by means of an explicit {\it
idealization} of any $\Dd$-triangulation of the triple. Moreover, the
explicit state sum expression $K(\Tt_N)$ of $K_N(W,L,\rho)$ tells us
that
\begin{equation}\label{conj1}
\lim_{N\to \infty} (2\pi/N^2) \log \left(|K_N(W,L,\rho)|\right) = 
G(\Tt_{\Ii})\ ,
\end{equation}
\noindent where $|\ |$ denotes the modulus of a complex number, and
$G$ basically depends on the geometry of the ideal tetrahedra of
$\Tt_{\Ii}$, and on the portion of the 1-skeleton of the triangulation
of $W$ which triangulates the link $L$. As $\Tt_{\Ii}$ is arbitrary,
roughly speaking again, $G$ may be considered as a function
$G(\mathfrak{c}_{\Ii}(W,L,\rho))$ of the $\Ii$-scissors congruence
class.

\noindent Following \cite{N1,N2}, there exists a refined version
$\widehat{\Pp(\Ii)}$ of the classical (pre)-Bloch group such that, by
using hyperbolic ideal triangulations of a non-compact complete and
finite volume hyperbolic $3$-manifold $M$ one can define a scissors
congruence class $\widehat{\beta}(M)\in \widehat{\Pp(\Ii)}$. Moreover,
one has
\begin{equation}\label{formrog}
R(\widehat{\beta}(M))= i(Vol (M) + iCS(M))\quad {\rm
mod}\left((\pi^2/2)\ \mz\right)\ ,
\end{equation} 
where $CS$ is the Chern-Simons invariant and $R:\widehat{\Pp(\Ii)} \to
\mc/\left((\pi^2/2)\ \mz\right)$ is a natural lift of the classical
Rogers dilogarithm on $\widehat{\Pp (\Ii)}$.  

\noindent
Starting from any $\Tt_{\Ii}$ as above, one can also define a refined
$\Ii$-class $\hat{\mathfrak{c}}_{\Ii}(W,L,\rho)$ and a 
{\it dilogarithmic invariant}
$$ R(W,L,\rho):= R(\hat{\mathfrak{c}}_{\Ii}(W,L,\rho)) \quad {\rm
mod}\left((\pi^2/2)\ \mz\right)\  .$$ 
\noindent There are strong structural relations between the classes 
$\mathfrak{c}_{\Ii}(W,L,\rho)$, $\hat{\mathfrak{c}}_{\Ii}(W,L,\rho)$
and $\widehat{\beta}(M)$. 
These relations and the actual asymptotic expansion of the cyclic $6j$-symbols 
support a formulation for triples $(W,L,\rho)$ of a so-called 
{\it Volume Conjecture}, which predicts, in particular, that
$$G(\mathfrak{c}_{\Ii}(W,L,\rho)) =
{\rm Im}\ R(\hat{\mathfrak{c}}_{\Ii}(W,L,\rho))\ . $$
\noindent In section \ref{vol} we discuss an extension of this
conjecture, involving the whole $K_N$ (not only $|K_N|$) and $R$. We
stress that the transition from $\Tt$ to $\Tt_{\Ii}$ is explicit and 
geometric and does not involve any ``optimistic'' computation. 
On the other hand the
actual identification of $G(\mathfrak{c}_{\Ii}(W,L,\rho))$ with ${\rm
Im}\ R(\hat{\mathfrak{c}}_{\Ii}(W,L,\rho))$ still sets serious analytic
problems.

\medskip
{\bf Acknowledgements}\qua We thank the referee for having forced us to
correct an early wrong formulation of the Complex Volume Conjecture
\ref{congecomplex}.

\section{$\Dd$-tetrahedra}\label{dtet}

We fix a {\it base} tetrahedron $\Delta$ embedded in $\mr^3$, with the
natural cell-decom\-posi\-tion. We orient $\mr^3$ by stipulating that
the standard basis is positive, and $\Delta$ is oriented in accordance
with it.  We consider $\Delta$ up to orientation-preserving cellular
self-homeomorphisms which induce the identity on the set of vertices.

\begin{figure}[ht!]
\begin{center}
\includegraphics[width=10cm]{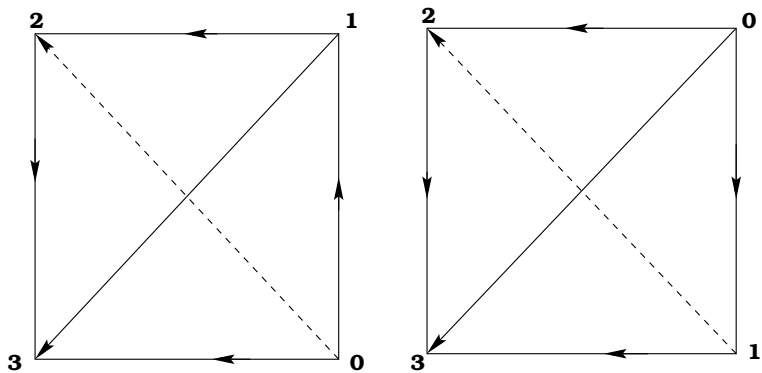}
\caption{\label{bstar} $+(\Delta,b)$ and $- (\Delta,b)$ }
\end{center}
\end{figure}

\noindent A $\Dd$-{decoration} on  $\Delta$ is a triple $((b,*),z,c)$ where:

\noindent (1)\qua $b$ is a {\it branching} on $\Delta$, that is a system of
edge-orientations such that no $2$-face inherits a coherent orientation on
its boundary.  It turns out that exactly one vertex is a source and
one vertex is a sink.  Every simplex of the natural triangulation of
the boundary of $\Delta$ inherits a ``name'' from the branching $b$:
if one denotes by $\Vv (\Delta)$ the set of vertices, these are named
by the natural ordering map
$$\{0,1,2,3\}\to \Vv (\Delta) \ \ \ \ \ \ \ i\to v_i$$ 
\noindent such that $[v_i,v_j]$ is an oriented edge of $(\Delta,b)$ if and
only if $j>i$.  Any other $j$-simplex is named by means of the names of
its vertices. The $2$-faces can be equivalently named in terms of
the opposite vertices.  One can select the ordered
triple of oriented edges
$$(e_0=[v_0,v_1],\ e_1=[v_1,v_2],\ e_2= -[v_2,v_0]=[v_0,v_2] )$$ 
\noindent which are the edges of the face opposite to the vertex $v_3$.  For every edge $e$ of
$\Delta$, one denotes by $e'$ the {\it opposite} edge. A branching
induces an orientation on $\Delta$ defined by the basis $(e_0',\ e_1',\ e_2')$, considered as
an ordered triple of vectors at $v_3$.  If this branching-orientation
agrees with the fixed orientation of $\Delta$ the branching is
said {\it positive}, and {\it negative} otherwise - see Figure
\ref{bstar}. We will write $*(\Delta,b)$, $*\in \{+,-\}$, to encode
the sign of the $b$-orientation of $\Delta$. We stipulate that $\bar{*}=-*$ with the
usual ``sign rule''. In a similar way, every simplex in $\partial
\Delta$ inherits an orientation from $b$.

\noindent (2)\qua $z$ is a {\it full} $B$-valued $1$-cocycle on
$(\Delta,b)$, where $B$ is the Borel subgroup of upper triangular
matrices of $SL(2,\mc)$. To specify $z$ one uses the edge-orientation
given by the branching $b$. For each oriented edge $e$, $z(e)$ is an
upper triangular matrix having $(t(e),1/t(e))$ on the diagonal and
upper diagonal entry equal to $x(e)$. So one can identify $z(e)$ with
$(t(e),x(e))$. ``Full'' means that $x(e)\neq 0$ for every $e$.  

\noindent (3)\qua $c$ is an {\it integral charge} on $\Delta$. Let $\Ee (\Delta)$
denote the set of edges of $\Delta$; then $ c: \Ee (\Delta) \to \mz $ is such that
\begin{itemize}
\item for every $e\in \Ee (\Delta)$, $c(e)=c(e')$;
\item If $e_i$, $i=0,\ 1,\ 2$, are the edges of any $2$-face of $\Delta$ and $c_i = c(e_i)$, then
$$ c_0+c_1+c_2 = 1\ .$$
\end{itemize}
It is useful (and pertinent - see \S \ref{preB}) to look at $c(e)\pi$ as a
dihedral angle; in this way the $c_i\pi$'s formally behave like the
dihedral angles of hyperbolic ideal tetrahedra.

\noindent Let us call $\Dd$ the set of all $\Dd$-tetrahedra
$*(\Delta,b,z,c)$. Denote by {\bf S}$_4$ the group of
permutations on $4$ elements. Changing the branching (i.e.\ permuting
the order of the vertices) induces a natural action
$p_{\Dd}$ of {\bf S}$_4$  on $\Dd$ with 
\begin{equation} \label{action}
p_{\Dd}(s,*(\Delta,b,z,c))= \epsilon(s)*(\Delta,s(b),s(z),s(c))\ ,
\end{equation}
\noindent where $\epsilon(s)$ is the signature of the permutation $s$, $s(b)$ is the new
branching obtained by permutation of the vertex ordering, and for every $e \in \Ee(\Delta)$ we have $s(z)(e)= z(e)$ iff the edge $e$ keeps the same orientation and $s(z)(e)= z(e)^{-1}$ otherwise, and $s(c)(e)=c(e)$.

\section{ The (pre)-Bloch-like group $\Pp (\Dd)$}\label{preB}
Let $\mz [\Dd]$ be the free $\mz$-module generated by the
$\Dd$-tetrahedra $*(\Delta,b,z,c)$; recall that we set
$\bar{*}(\Delta,b,z,c)= (-1)*(\Delta,b,z,c)$.

\noindent 
In this section we shall describe a notion of $2 \leftrightarrow 3$
$\Dd$-\emph{transit} between decorated triangulations such that every
instance of $2 \leftrightarrow 3$ $\Dd$-transit produces a natural {\it five term}
relation in $\mz [\Dd]$. By definition $ \Pp (\Dd) = \mz [\Dd]/T(\Dd)$
is the {\it (pre)-Bloch-like} $\Dd$-{\it group}, where $T(\Dd)$ is the
module generated by all these five term relations and by the relations
(\ref{action}).

\begin{figure}[ht!]
\begin{center}
\includegraphics[width=12cm]{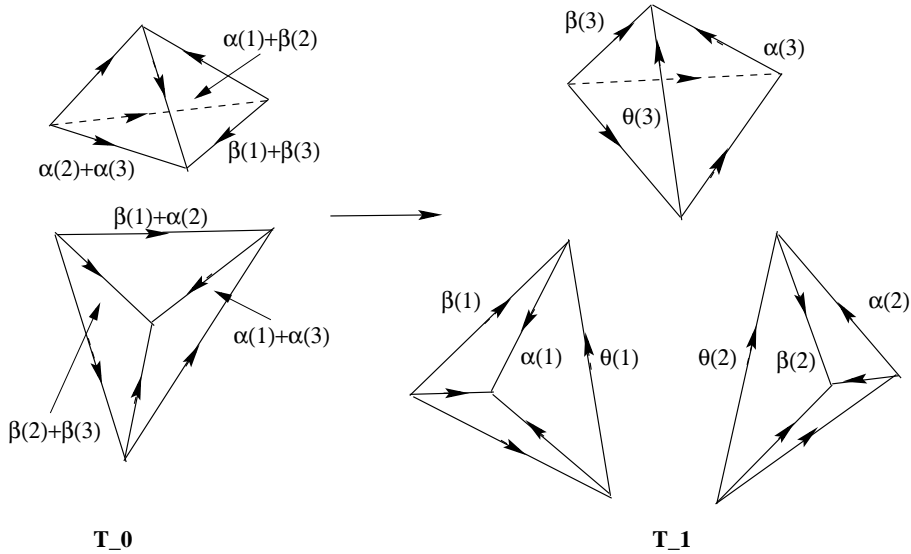}
\caption{\label{charget} An instance of branching and charge transit}
\end{center}
\end{figure}

\noindent 
The support of any $2 \leftrightarrow 3$ $\Dd$-transit is the usual $2 \leftrightarrow 3$ move
on $3$-dimensional triangulations. Let us specify how the decorations
transit. In Figure \ref{charget}, forgetting the charge for a while,
one can see an instance of {\it branching transit}. The branching
transits are also carefully analyzed in \cite{BP}, in terms of the
dual viewpoint of ``branched spines''.

\noindent 
A {\it charge transit} is branching independent. Consider a $2 \leftrightarrow 3$
move $T_0 \leftrightarrow T_1$ relating two triangulations of some $3$-manifold,
and suppose that the tetrahedra involved in the move are endowed with
integral charges. There are $9$ edges $E_1, \dots, E_9 \ $ which are
present in both $T_0$ and $T_1$ and a further edge $E_0$ which is
present only in $T_1$. Also, for $i=1,\dots,9 \ $, $E_i$ is an edge of
exactly one tedrahedron of $T_k$ iff it is an edge of
exactly two tetrahedra of $T_ {k+1}$ (where $k \in \mz/2\mz$), and
$E_0$ is an edge of exactly 3 tetrahedra of $T_1$.  Denote by
$\gamma_k (E_i)$, $i=1,\dots,9 \ $, the sum of the integral charges at
$E_i$ of the tetrahedra of $T_k$ which have $E_i$ among their edges
(by convention put $\gamma_0 (E_0)=0$). The integral charges on $T_0$ and $T_1$
define a $2 \leftrightarrow 3$ {\it charge transit} iff they coincide for the
tetrahedra not involved in the move, and the following relations are
satisfied:
$$   \gamma_0 (E_i) = \gamma_1 (E_i) \ ,\quad i=1,\dots,9\ .$$ 
These linear relations imply $\gamma_1 (E_0) = 2$, which together
with the second condition in
\S \ref{dtet} (3) for the integral charges on $T_1$ read:
$$\theta(1) + \theta(2) + \theta(3) = 2\quad ;\quad \theta(j) +
\alpha(j) + \beta(j) = 1 \ ,\quad j=1,2,3 \ .$$

These two relations show that a charge transit is coherent 
with the above ``dihedral angle'' interpretation of the integral charges.

\begin{figure}[ht!]
\begin{center}
\includegraphics[width=12cm]{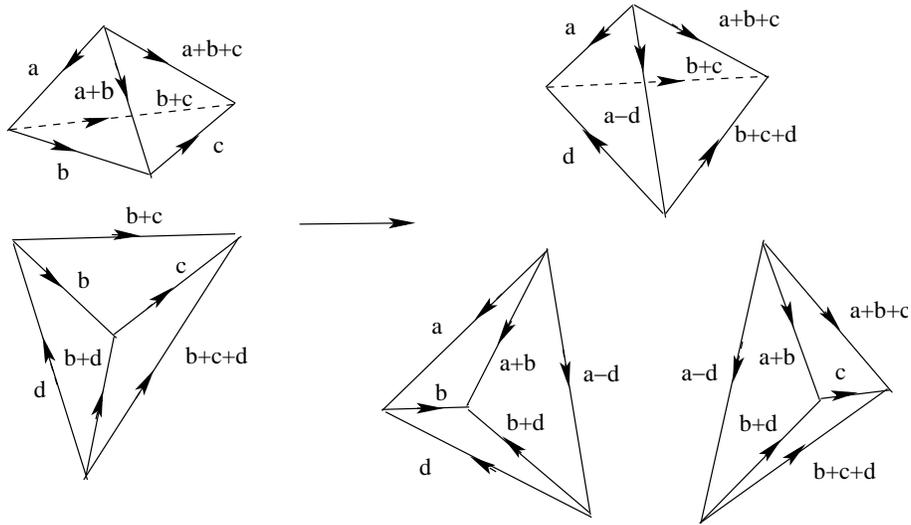}
\caption{\label{cyclet} An instance of cocycle transit}
\end{center}
\end{figure}

\noindent The last decoration component of a $\Dd$-decoration is the
cocycle $z$. Let $T_0$, $T_1$ be as above and $z_0 \in Z^1(T_0;
B)$. Then there is only one $1$-cocycle $z_1$ in $Z^1(T_1; B)$ such
that $z_0(E_i)=z_1(E_i)$ for $i=1,\dots, 9\ $. It defines the {\it
cocycle transit}, providing that both cocycles $z_0$ and $z_1$ are full. 
In Figure \ref{cyclet} one can see an instance of cocycle transit where,
for the sake of simplicity, we have used cocycles in $Z^1(T_k;\mc)\cong
Z^1(T_k;Par(B))$, where $Par(B)$ is the parabolic subgroup of $B$ of matrices
with unitary diagonal.

\noindent 
Any instance of $2 \to 3$ $\Dd$-transit induces a five term relation
in $\mz [\Dd]$; note that we have to take into account the
orientations (i.e.\ the signs) of the $5$ involved tetrahedra in these
relations.

\section{ The $\Dd$-scissors congruence class} \label{scissors}
Let $(W,L,\rho)$ be a triple formed by a smooth compact closed
oriented $3$-manifold $W$, a non-empty link $L$ in $W$ and a flat principal 
$B$-bundle $\rho$ on $W$. The triple is regarded up to
orientation-preserving diffeomorphisms of  $(W,L)$ which are flat bundle isomorphisms.

\noindent 
A {\it distinguished triangulation} $(T,H)$ of $(W,L)$ is a (singular)
triangulation $T$ of $W$ such that $L$ is realized by a {\it Hamiltonian}
sub-complex $H$ of $T$ that contains all the vertices. Such a
triangulation can be interpreted as a certain finite family
$\{\Delta_1,\dots, \Delta_k\}$ of copies of our base tetrahedron
endowed with a set of $2$-faces identification (pairing) rules.  A
{\it $\Dd$-decoration} $(b,z,c)$ on $(T,H)$ consists of a family of
$\Dd$-tetrahedra $\{*_i(\Delta_i,b_i,z_i,c_i)\}$, $i=1,\dots,k$, which
is compatible with the face pairings; this means that

\noindent 
(1)\qua the pairings respect the edge-orientations due to the branchings;
    thus they give us a global branching $b$ on $T$ (see
    \cite{BP});

\noindent 
(2)\qua identified edges have the same cocycle value, so that the $z_i$'s
    define a global $B$-valued full $1$-cocycle $z$ on $T$;

\noindent 
(3)\qua denote by $\textstyle \Ee = \coprod_i \Ee(\Delta_i)$ the set of
all edges of all $\Delta_i$'s; the system of integral charges $c_i$
can be considered as a global integral charge $c: \Ee \to \mz$.

\noindent Moreover, one imposes the
following  global constraints:

\noindent (4)\qua a branching $b_i$ is positive if and
only if the corresponding $b_i$-orientation on $\Delta_i$ agrees with 
the one of the manifold $W$;
\smallskip

\noindent (5)\qua the global cocycle $z$ on $T$ represents the flat bundle
$\rho$ on $W$: $\rho = [z]$;
\smallskip

\noindent (6)\qua denote by  $\Ee (T)$ the set of edges of $T$. 
There is a natural projection map $p: \Ee \to \Ee (T)$. Then one requires that
for every $s\in \Ee (T)\setminus \Ee (H)$
$$ \sum_{e \in p^{-1}(s)} c(e) = 2\ ,$$
and for every $s\in \Ee (H)$
$$ \sum_{e \in p^{-1}(s)} c(e) = 0\ .$$
\noindent (7)\qua to every $c$ which satisfies condition (6) one can
associate an element $[c]\in H^1(W;\mz/2\mz)$; one finally imposes that
$[c]=0$. 

\smallskip

\noindent Note that, by property (6), any integral charge on $(T,H)$ actually
encodes $H$. We say that that a triangulation $T$ of $W$ is {\it fullable}
if it carries a full $1$-cocycle representing the trivial flat
$B$-bundle. This is equivalent to the fact that every edge of $T$ has
two distinct vertices. If $T$ is fullable then for any flat $B$-bundle
$\rho$ it carries a full $1$-cocycle representing $\rho$. Finally
$\Tt=(T,H,(b,z,c))$ is called a $\Dd$-{\it triangulation} of $(W,L,\rho)$ if
$T$ is fullable and $z$ is full.

\begin{theo}\label{exist} 
For every triple $(W,L,\rho)$ there exist $\Dd$-triangulations.
\end{theo}

\noindent To any $\Dd$-triangulation $\Tt$
one can associate the formal sum $\mathfrak{c}_{\Dd}(\Tt) \in \mz
[\Dd]$ of its decorated tetrahedra (all the coefficients being equal
to $1$).

\begin{theo}\label{Dclass} 
The equivalence class of $\mathfrak{c}_{\Dd}(\Tt)$ in $\Pp(\Dd)$ does
not depend on the choice of $\Tt$. Thus it defines an element
$\mathfrak{c}_{\Dd}(W,L,\rho)\in \Pp(\Dd)$, which is called the {\rm
$\Dd$-scissors congruence class} (or {\rm $\Dd$-class}) of the triple
$(W,L,\rho)$.
\end{theo}

\noindent 
The proof of this theorem is similar to the proof of the invariance of
the QHI state sums (see the next section). However, the proof requires
more triangulation moves (such as the so called {\it bubble} move)
than the only $2 \leftrightarrow 3$ move. A remarkable fact is that these further
moves do not introduce new independent algebraic relations in $\mz
[\Dd]$.

\section{Quantum hyperbolic state sum invariants}\label{QHI}

\noindent 
Let $(W,L,\rho)$ be as in \S \ref{scissors} and fix any 
$\Dd$-triangulation $\Tt=(T,H,(b,z,c))$ of $(W,L,\rho)$.  Let $N>1$ be
an odd integer. Fix a determination of the $N$th-root which holds for
all the matrix entries $t(e)$ and $x(e)$ of $z(e)$, for all the edges
$e$ of $T$.  The {\it reduction mod$(N)$} $\Tt_N$ of $\Tt$ consists in
changing the decoration of each edge $e$ of $T$ as follows:
\begin{itemize}
\item $(a(e)=t(e)^{1/N},\ y(e)=x(e)^{1/N})$ instead of $z(e)=(t(e),x(e))$; 
\item $c_N(e)=c(e)/2$ mod$(N)$ instead of $c(e)$ (it makes sense
because $N$ is odd).
\end{itemize}
\noindent Let us interpret this new decoration. All the details and
justifications of the explicit formulas given below can be found in
\cite{B}, see also the appendix of \cite{BB1}. Consider a quantum
Borel subalgebra $\Ww_N$ of $U_q(sl(2,\mc))$, specialized at the root
of unit $\omega_N = \exp (2\pi i/N)$. Each $(a(e),y(e))$ describes an
irreducible $N$-dimensional {\it cyclic representation} $r_N(e)$ of
this algebra.

\noindent Call $\Ff (T)$ the set of $2$-faces of $T$.  A {\it
$N$-state} of $T$ is a function $\alpha : \Ff (T) \to \{0,1,\dots,\
N-1\}$ (in fact one often identifies $\{0,1,\dots,\ N-1\}$ and
$\mz/N\mz$). The state $\alpha$ can be considered as a family of
functions $\alpha_i : \Ff(\Delta_i)\to \{0,1,\dots,\ N-1\}$ which are
compatible with the face pairings. 

\noindent Consider on each branched tetrahedron
$(\Delta_i,b_i)$ of $\Tt$ the ordered triple of oriented edges
$(e_0=[v_0,v_1],\ e_1=[v_1,v_2],\ e_2= -[v_2,v_0])$ which are the opposite
edges to the vertex $v_3$.  The cocycle property of $z$ and the fullness
assumption (this is crucial at this point, due to the algebraic
structure of the cyclic representations of $\Ww_N$) imply that
$r_N(e_0)\otimes r_N(e_1)$ coincides up to isomorphism with the direct
sum of $N$ copies of $r_N(e_2)$.  This set of data allows one to
associate to every $*(\Delta_i,b_i,r_{N,i},\alpha_i)$ a $6$j-symbol
$R(*(\Delta_i,b_i,r_{N,i},\alpha_i))\in \mc$, that is a matrix element
of a suitable ``intertwiner'' operator.  The reduced charge $c_N$ is
used to slightly modify this operator in order to get its (partial)
invariance up to branching changes. In this way one
gets the (partially) symmetrized {\it $c$-$6$j-symbols}
$T(*(\Delta_i,b_i,r_{N,i},c_i,\alpha_i)) \in \mc$.  We are now ready
to define the state sums $H(\Tt_N)$, which are {\it weighted} traces. 
Denote by $V$ the number of vertices of $T$. Set
$$ \Psi (\Tt_N)= \sum_\alpha \prod_i
T(*(\Delta_i,b_i,r_{N,i},c_i,\alpha_i))$$
\begin{equation} \label{inv}
H(\Tt_N) = \Psi (\Tt_N) \ N^{-V} \prod_{e\in \Ee (T)\setminus \Ee
(H)} x(e)^{(N-1)/N}\ .
\end{equation}

\begin{theo}\label{statesum}  
Up to multiplication with $N$th-roots of unity, the scalar $H(\Tt_N)$
does not depend on the choice of $\Tt$. Hence $K(\Tt_N) := H(\Tt_N)^N$
defines an invariant $K_N(W,L,\rho)$ of the triple $(W,L,\rho)$.
\end{theo}

\begin{proposi}\label{proprieta}  Let $\rho = [z]=[(t,x)]$.
\begin{itemize}
\item For every $\lambda \ne 0$ set $\lambda
\rho = [(t,\lambda x)]$. Then $K_N(W,L,\rho) = K_N(W,L,\lambda \rho)\ .$ 
\item Denote by $z^*$ the complex conjugate cocycle of $z$ and by $\rho ^*=[z^*]$
the corresponding flat $B$-bundle. Let $-W$ be $W$ endowed with the opposite
orientation. Then $K_N(-W,L,\rho^*)= K_N(W,L,\rho)^*\ .$
\end{itemize}
\end{proposi}

\noindent Kashaev proposed in \cite{K1} a conjectural purely
topological invariant $K_N(W,L)$, which should have been expressed by
a state sum as in  (\ref{inv}). In fact, one eventually
recognizes such a $K_N(W,L)$ as the special case of our
$K_N(W,L,\rho)$ when $\rho$ is the {\it trivial} flat $B$-bundle on
$W$ (although in \cite{K1} there are neither flat bundles, nor
geometric interpretations nor existence results of all the datas).

\noindent The algebraic properties of the c-$6j$-symbols (the
``pentagon relation'' and so on) ensure the invariance of $K(\Tt_N)$
up to $\Dd$-transits supported by certain ``bare'' triangulation
moves ($2 \leftrightarrow 3$ moves as defined in \S \ref{preB}, ``bubble
moves'',$\dots$). Rephrased in our setup, this was the main
achievement of \cite{K1}. However, one cannot deduce the complete
invariance of $K(\Tt_N)$ solely from this $\Dd$-transit invariance because
it is difficult to connect by $\Dd$-transits two
given $\Dd$-triangulations of a triple
$(W,L,\rho)$.  For example, ``negative'' $3\to 2$ moves are in general
not ``brancheable'', and full cocycles do not transit, in general, to
full cocycles. Also the charge invariance relies deeply on the fact that
the set of integral charges of any fixed distinguished triangulation
is an integral lattice.

\noindent 
The nature of the ambiguity of $H(\Tt_N)$, up to $N$th-roots of unity,
is not yet clear to the authors. The problem is due to the
symmetrization procedure, which turns $6j$-symbols into
c-$6j$-symbols. Indeed, the $6j$-symbols have a very subtle behaviour
w.r.t. branchings: only the pentagon relations corresponding to a
non-trivial proper subset of branched $2 \leftrightarrow 3$ moves are valid. 

\noindent 
As the value of $K(\Tt_N)$ does not depend on the choice of $\Tt$ one
would like to consider $K_N(W,L,\rho)$ as a function of the
$\Dd$-class $\mathfrak{c}_{\Dd}(W,L,\rho)$. This is not completely
correct because the face pairings between the $\Dd$-tetrahedra of
$\Tt$ are not encoded in the representatives $\mathfrak{c}_{\Dd}(\Tt)$
of $\mathfrak{c}_{\Dd}(W,L,\rho)$. Moreover, the states as well as the
non-$\Psi(\Tt_N)$ factors in the right-hand side of (\ref{inv}) depend
on the face pairings. This is a technical point which can be overcome
as follows, by looking at $K(\Tt_N)$ as a function of a formal sum of
``augmented'' $\Dd$-tetrahedra.

\noindent A $\widetilde{\Dd}$-{\it tetrahedron} is of the form
$*(\Delta,b,z,c,v^0,v^1,v^2)$ where $v^0,v^1,v^2$ are $\mn$-valued
functions defined on $\Vv (\Delta),\ \Ee (\Delta),\ \Ff (\Delta)$
respectively. Let $\Gamma \in \mz [\widetilde{\Dd}]$ be a formal sum of
terms with coefficients all equal to $1$. For every $\sigma : \mn \to
\mz/N\mz$ set $\alpha_i(\sigma)=\sigma \circ v_i^2$. Say that $\sigma$ and
$\sigma'$ are ``identified modulo $\Gamma$'' if $\alpha_i =
\alpha_i'$. Put
$$\Phi = \prod_{w\in \Vv (\Delta_i)} N^{-1/v_i^0(w)}\quad , \quad \Omega = \prod_ {e\in \Ee (\Delta_i)} x(e)^{(N-1)/v_i^1(e)N}\ $$
$$\widetilde{T}(*(\Delta_i,b_i,r_{N,i},c_i,\alpha_i(\sigma),v_i^0,v_i^1))= 
T(*(\Delta_i,b_i,r_{N,i},c_i,\alpha_i(\sigma)) \ \Phi \ \Omega \ .$$

\noindent Finally set
$$ H(\Gamma) = \sum_{\sigma \in (\mz/N\mz)^{\mn}/\Gamma} \prod_i
\widetilde{T} (*(\Delta_i, b_i, r_{N,i}, c_i, \alpha_i(\sigma),v_i^0,
v_i^1))\ .$$
\noindent Let $\Tt$ be a $\Dd$-triangulation. For every vertex $w$ of $\Delta_i$ set $v_i^0(w) =
|p_0^{-1}p_0(w)|\in \mn$, where $p_0$ is the identification map in $T$ of the vertices of the $\Delta_i$'s. For every edge $e$ of $\Delta_i$ set $v_i^1(e) = |p^{-1}p(e)|\in \mn$, where $p: \Ee \to \Ee(T)$ is as in \S \ref{scissors} (6). Finally order $\Ff (T)$ via an arbitrary $\mn$-valued map which is compatible (for the face pairings) with $v^2$. In this way we get a
$\widetilde{\Dd}$-triangulation $\widetilde{\Tt}$ with the associated
$\mathfrak{c}_{\widetilde{\Dd}}(\widetilde{\Tt}) \in \mz [\widetilde{\Dd}]$.  It is clear that
$$K_N(W,L,\rho)=H(\mathfrak{c}_{\widetilde{\Dd}}(\widetilde{\Tt}))^N\ .$$  
\noindent Following the lines of the above construction,
it is not hard to define the notions of $\widetilde{\Dd}$-transit, (pre)-Bloch-like group
$\Pp (\widetilde{\Dd})$ and $\widetilde{\Dd}$-class
$\mathfrak{c}_{\widetilde{\Dd}}(W,L,\rho)\in \Pp (\widetilde{\Dd})$, so that one can
eventually define a map $K_N: \Pp (\widetilde{\Dd}) \to \mc$ such
that
$$ K_N(W,L,\rho) = K_N(\mathfrak{c}_{\widetilde{\Dd}}(W,L,\rho))\ .$$ 
\noindent Moreover, there is a ``forgetting map'' $\Pp (\widetilde{\Dd}) \rightarrow \Pp (\Dd)$ sending
$\mathfrak{c}_{\widetilde{\Dd}}(W,L,\rho)$ to $\mathfrak{c}_{\Dd}(W,L,\rho)$. With this 
technical precision in mind, we shall say anyway that, roughly speaking, $K_N(W,L,\rho)$ depends on $\mathfrak{c}_{\Dd}(W,L,\rho)$.

\section{The (pre)-Bloch like group $\Pp (\Ii)$}  \label{preid}

\noindent Let us go back to our base tetrahedron $\Delta$. An $\Ii$-{\it
decoration} on $\Delta$ is a triple $((b,*),w,c)$ where $(b,*)$ and
$c$ are as in the $\Dd$-decorations and $w: \Ee (\Delta) \to \mc \setminus \{0,1\}$ is such that :
\begin{itemize}
\item for every $e\in \Ee (\Delta)$, if $e'$ is opposite to $e$, then $w(e)=w(e')$;
\item if $e_0,\ e_1,\ e_2$ are the edges of the face opposite
to $v_3$ and $ w_i = w(e_i)$, then $ w_0w_1w_2 = -1\ \ \ {\rm and}\ \ \ w_0w_1 - w_1 = -1 $.
\end{itemize}

\noindent 
Clearly, $w'=w_1 = (1-w_0)^{-1}$, $w''=w_2 = (1-w_1)^{-1}$ and $w=w_0
= (1-w_2)^{-1}$, so that $(w,w',w'')$ is the {\it modular triple} of
an oriented hyperbolic ideal tetrahedron $\bar{\Delta}$ in oriented
$\mh^3$. If ${\rm Im}(w)>0$ then $\bar{\Delta}$ is oriented
positively, and negatively otherwise. However, we do not require that
the $b$-orientation of $*(\Delta,b,w,c)$ coincides with the sign of
${\rm Im}(w)$.

\begin{figure}[ht!]
\begin{center}
\includegraphics[width=12cm]{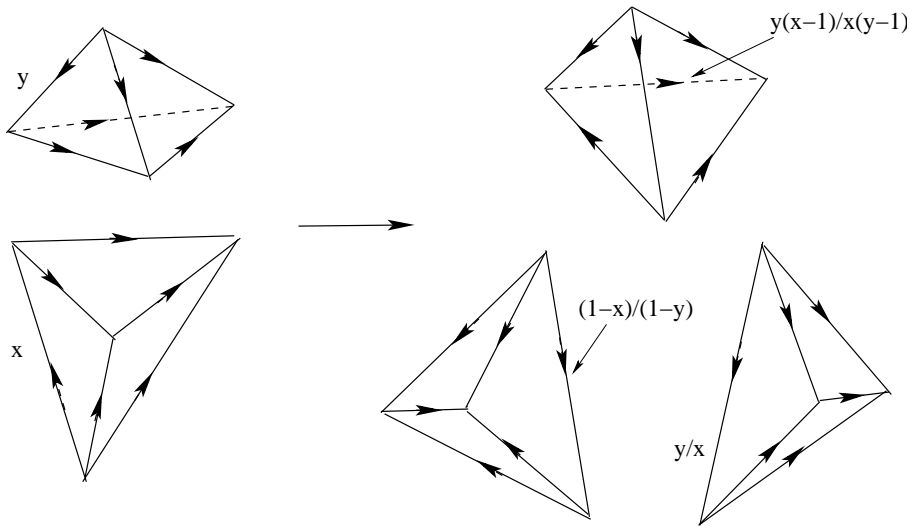}
\caption{\label{idealt} An instance of ideal transit}
\end{center}
\end{figure}

\noindent 
Figure \ref{idealt} shows an instance of $2 \leftrightarrow 3$
$\Ii$-{\it transit}. Only the first members of the modular triples are
indicated. Note that we are assuming that both on $T_0$ and $T_1$ we
actually have modular triples; this means that the $\Ii$-transit is possible
only if $x \ne y$. This fact is strictly related to the fullness requirement for $\Dd$-transits (see \S \ref{idealization}).  
Each instance of $\Ii$-transit operates on the branching
and on the integral charges like a $\Dd$-transit. Recall that a
modular triple determines and is determined by the dihedral angles at
the edges of the corresponding ideal tetrahedra. Then, in terms of
dihedral angles, modular triple transits are formally defined like
integral charge transits. In another (equivalent) way, in a $2
\leftrightarrow 3$ $\Ii$-transit that splits $*(\Delta,b,w)$ into
$*_1(\Delta_1,b_1,w_1)$ and $*_2(\Delta_2,b_2,w_2)$, we have $w(e)
^{*}=w_1(e)^{*_1}w_2(e)^{*_2}$ (recall that each of $*,*_1$ and $*_2$
equals $\pm 1$).

\noindent Denote by $\Ii$ the set of all $\Ii$-tetrahedra
$*(\Delta,b,w,c)$. Any instance of $2 \leftrightarrow 3$ $\Ii$-transit
produces a five term relation in the free $\mz$-module $\mz
[\Ii]$. Also, there is a natural action $p_{\Ii}$ of {\bf S}$_4$ on
$\Ii$ which acts as $p_{\Dd}$ in (\ref{action}) on $b$, $*$ and $c$;
moreover $s(w)(e)=w(e)^{\epsilon (s)}$. One defines the {\it
(pre)-Bloch-like group} $ \Pp (\Ii)$ as the quotient of $\mz [\Ii]$ by
the module generated by all the $2 \leftrightarrow  3$ $\Ii$-transit five terms
relations and by the relations induced by $p_{\Ii}$.

\section{Idealization of $\Dd$-tetrahedra}\label{idealization}

\noindent 
Let $*(\Delta,b,z,c)\in \Dd$, and $(e_0=[v_0,v_1],\ e_1=[v_1,v_2],\
e_2= -[v_2,v_0])$ be as in \S \ref{dtet}. The cocycle property of
$z=(t,x)$ implies that
\begin{equation} \label{trucmuch}
  x(e_0)x(e_0')+ x(e_1)x(e_1')+ (-x(e_2) x(e_2')) := p_0 + p_1 + p_2=0\ .
\end{equation}
\noindent Define (indices mod($\mz_3$)):
$$ F(*(\Delta,b,z,c)) = *(\Delta,b,w(z),c), \quad w_i= -
p_{i+1}/p_{i+2}\ .$$ 
\noindent It is readily seen that
$F(*(\Delta,b,z,c))$ belongs to $\Ii$; we call $F$ the {\it
idealization map}. In fact, consider an oriented hyperbolic ideal
tetrahedron $\bar{\Delta}$ with vertices $v_j \in S_{\infty}^2 =
\partial \bar{\mathbb{H}}^3$ and modular triple $(w_0,w_1,w_2)$. One
can assume that every $v_j \in \mc \subset \mc \cup \{\infty \}=
\partial \mh^3$.  If one looks at the $v_j$'s as defining a $0$-cochain $u$ on $\bar{\Delta}$, then the usual {\it cross-ratio} expressions 
$$w_0 = [v_0:v_1:v_2:v_3]=
\frac{(v_2-v_1)(v_3-v_0)}{(v_2-v_0)(v_3-v_1)}\ ,$$
\noindent etc., are compatible with the definition of $F$ by using the 1-cocycle given by $\delta(u)$. Note also that the idealization only depends on the ``projective''
class of the full cocycle, similarly to the behaviour of QHI in
Proposition \ref{proprieta}. 
Clearly, the map $F$ is onto; it extends linearly to $F: \mz
[\Dd] \to \mz [\Ii]$. Remarkably one has
 
\begin{proposi}\label{parB} The map $F$ induces a well-defined surjective 
homomorphism
$ \widehat{F}: \Pp (\Dd)\to \Pp (\Ii)$.
\end{proposi}

\section{Hyperbolic-like structures -- 
A volume conjecture}\label{vol} 

\noindent 
A fundamental problem in QHI theory is to understand the asymptotic
behaviour of the state sum invariants $K_N(W,L,\rho)$ when $N\to
\infty$. Accordingly with the considerations of \S \ref{QHI}, this
should depend on the $\Dd$-class of $(W,L,\rho)$. Note that one has
\begin{equation}\label{asymptot} \lim_{N\to 
\infty} (2\pi/N^2) \log (|K_N(W,L,\rho)|) = \lim_{N\to \infty}
(2\pi/N) \log(|\Psi (\Tt_N)|)
\end{equation}
\noindent This limit is finite and does not vanish iff
$|K_N(W,L,\rho)|$ grows exponentially w.r.t. $N^2$. This is
corroborated by the computation of the asymptotic behaviour of the
c-$6j$-symbols; their exponential growth rate
involves classical dilogarithm functions.  Moreover, the
explicit expression of the c-$6j$-symbols shows that
\begin{equation}\label{property}
T(*(\Delta_i,b_i,r_{N,i},c_i,\alpha_i)) =
T(F(*(\Delta_i,b_i,r_{N,i},c_i,\alpha_i)))\ .
\end{equation}
If $\Tt$ is any $\Dd$-triangulation, then $F(\Tt)$ is an
$\Ii$-triangulation. Using it, one can define a conjugacy class of
holonomy representations $\hat{\rho}_{\Tt}: \pi_1(W) \rightarrow
PSL(2,\mc)$ equipped with piecewise straight equivariant
maps from the universal covering of $W$ to $\bar{\mh}^3$. Moreover,  set
$\mathfrak{c}_{\Ii}(F(\Tt))=\widehat{F}(\mathfrak{c}_{\Dd}(\Tt)) \in
\Pp (\Ii)$. One has

\begin{theo} \label{hlt} 
Both $\hat{\rho}_{\Tt}$ and $\mathfrak{c}_{\Ii}(F(\Tt))$ do not depend
on $\Tt$. Hence they define respectively an {\rm hyperbolic-like
structure} $\hat{\rho}$ on $W$ which only depends on 
$(W,\rho)$, and an $\Ii$-{\rm scissors congruence class} 
(shortly {\rm $\Ii$-class}) 
$\mathfrak{c}_{\Ii}(W,L,\rho)$. 
\end{theo}

\noindent The proof of Theorem \ref{hlt} essentially follows from
Theorem \ref{Dclass} and \ref{parB}. It also uses the following
remarkable geometric feature. Let $\{*(\Delta_i,b_i,w_i,c_i)\}$ be the
family of $\Ii$-tetrahedra of $F(\Tt)$. The union of $w_i$'s can be
considered as a map $w$ defined on the set $\Ee$ of edges of all
$\Delta_i$'s. Let $p:\Ee \to \Ee (T)$ be as in \S \ref{scissors}
(6). Then, for every $s\in \Ee (T)$, one has $\textstyle \prod_{ e\in
p^{-1}(s)} w(e)^*=1$.

\noindent 
If $\Tt$ is a $\Dd$-triangulation, using (\ref{property}) one
can rewrite (\ref{asymptot}) as
$$ \lim_{N\to \infty} (2\pi/N^2) \log (|K_N(W,L,\rho)|) = G(F(\Tt )) \
.$$ Here the right-hand side is a function of the hyperbolic ideal
tetrahedra of $F(\Tt)$ and of $H$, which is encoded by the charge in
$\Tt$. As $G(F(\Tt )) $ does not depend on the choice of $\Tt$, roughly speaking
(see the discussion at the end of \S \ref{QHI}), $G$ is a function of
the $\Ii$-class of $(W,L,\rho)$.

\noindent
Starting from any idealized triangulation $F(\Tt)$ of $(W,L,\rho)$ one
can also define a refined $\Ii$-class
$\hat{\mathfrak{c}}_{\Ii}(W,L,\rho)$.  This class can be represented
by a certain decorated triangulation $F'(\Tt)$ which differs from
$F(\Tt)$ by adding a ``combinatorial flattening'' (see  \cite{N2}) in 
the decoration. The trace of the latter on each
tetrahedron of $T$ behaves as a signed charge, and it satisfies 
global constraints similar to conditions (6) and (7) in Section \ref{scissors},
but depending also on the moduli.

\noindent Following the comments in the Introduction, one then defines a 
{\it dilogarithmic invariant} as
$$ R(W,L,\rho):= R(\hat{\mathfrak{c}}_{\Ii}(W,L,\rho)) \quad {\rm
mod}\left((\pi^2/2)\ \mz\right)\  .$$
A first formulation of the Volume Conjecture is:

\begin{conge}[Real Volume Conjecture]\label{congmodule} 
For any $(W,L,\rho)$ we have
$$ \lim_{N\to \infty} (2\pi/N^2) \log \left(| K_N(W,L,\rho)|\right) = 
{\rm Im}\ R(\hat{\mathfrak{c}}_{\Ii}(W,L,\rho))\ .$$
\end {conge}
\noindent Conjecture \ref{congmodule} is in {\it formal} agreement
with the current Volume Conjecture for hyperbolic knots in $S^3$ based
on the coloured Jones polynomial $J_N$ (\cite{Kbis,K2},\cite{MM,Y}),
because it is commonly accepted that $(J_N)^N$ is an instance of
$K_N$. However, in our opinion, this has not yet been proved anywhere
and the question of the relationship between $K_N$ and $J_N$ needs
further investigation (see \cite{BB2}).

\noindent A rough idea to extend the above conjecture for
$K_N(W,L,\rho)$ (not only for its modulus) is to formally pass to an
exponential version of Conjecture \ref{congmodule}, and replace ${\rm
Im}\ R(\hat{\mathfrak{c}}_{\Ii}(W,L,\rho))$ with
$R(\hat{\mathfrak{c}}_{\Ii}(W,L,\rho))$. But
$R(\hat{\mathfrak{c}}_{\Ii}(W,L,\rho))$ is only determined mod($(\pi^2/2) \
\mz$). One avoids this ambiguity as follows:

\begin{conge}[Complex\ \ Volume\ \ Conjecture]\label{congecomplex} 
There exist invariants\break $C(W,L,\rho) \in \mc \ {\rm mod}((\pi^2/2)\ \mz)$ and $D=D(W,L,\rho)\in \mc^*$ such that for any branches
$R$ and $C$ of $R(\hat{\mathfrak{c}}_{\Ii}(W,L,\rho))$ and $C(W,L,\rho)$ respectively, we have
$$ K_N(W,L,\rho)^8 = \left[ \exp\left( \frac{C+ N R}{2i\pi}\right) \right]^{8N} \left(D + \mathcal{O}(\frac{1}{N})\right)\ .$$  
\end{conge}
\noindent Conjecture \ref{congecomplex} says at first that
$K_N(W,L,\rho)^8$ has an exponential growth rate. Assuming it, the
fact that $\exp (4C/i\pi)$, $\exp(4R/i\pi)$ and $D$ are
well-determined invariants of $(W,L,\rho)$ follows from the invariance
of $ K_N(W,L,\rho)$ and the uniqueness of the coefficients of
asymptotic expansions (of Poincar\'e type). Moreover Conjecture
\ref{congecomplex} predicts that $R$ is a branch of
$R(\hat{\mathfrak{c}}_{\Ii}(W,L,\rho))$. 

\noindent At present, the nature of
$C(W,L,\rho)$ and $D(W,L,\rho)$ is somewhat mysterious. There are no reasons to
expect that, for instance, $C(W,L,\rho)=0$ or $D(W,L,\rho)=1$; in fact, their
value could be related to the hard problem of finding an appropriate
contour of integration in the stationary phase approach to the
evaluation of the asymptotic behaviour of $K_N(W,L,\rho)$.
 
\noindent A natural complement to the conjecture is the problem of  
understanding the geometric meaning of the dilogarithmic invariant of a triple
$(W,L,\rho)$.

\Addresses\recd

\end{document}